\documentclass[12pt]{article}
\usepackage{amssymb,amsthm,amsmath}

\theoremstyle{plain}
\newtheorem{theorem}{Theorem}

\title{Spectral portraits and the resolvent growth of a model problem
associated with the Orr--Sommerfeld equation}
\author{Shkalikov~A.~A.\thanks{Research supported by RFFI grant No~01-01-00691
and by the grant "The Universities of Russia"}\\
Department of Mechanics and Mathematics\\
Moscow State University\\
\texttt{ashkalikov@yahoo.com}}
\date{}
\begin{document}
\maketitle
The Orr--Sommerfeld equation is obtained by linearization of the Navier--Stokes
system in the infinite spacial layer \((x,\xi,\eta)\in\mathbb R^3\), where
\(|x|\leqslant 1\) and \((\xi,\eta)\in\mathbb R^2\). This equation has the
form
\begin{equation}\label{eq1}
	(D^2-\alpha^2)y-i\alpha R\left[q(x)(D^2-\alpha^2)-q''(x)\right]y=
	i\alpha R\lambda(D^2-\alpha^2)y.
\end{equation}
Here \(D=d/dx\), \(\alpha\) is the wave number (\(\alpha\neq 0\)) appearing
after
the separation of the variables \((\xi,\eta)\in\mathbb R^2\), \(R\) is the
Reynolds number, \(q(x)\) is the stationary unperturbed velocity profile and
\(\lambda\) is the spectral parameter. Usually the boundary conditions
\begin{equation}\label{eq2}
	y(\pm 1)=y'(\pm 1)=0
\end{equation}
are associated with equation~\eqref{eq1}.

Our main goal is to understand the spectrum behaviour of
problem~\eqref{eq1},~\eqref{eq2} as the Reynolds number \(R\) tends to the
infinity (this means that the liquid becomes to be ideal). For these
purposes we associate with~\eqref{eq1},~\eqref{eq2} a problem of the form
\begin{gather}\label{eq3}
	-i\varepsilon^2 z''+q(x)z=\lambda z\\ \label{eq4}
	z(-1)=z(1)=0,
\end{gather}
assuming that \(q(x)\) is the same as in~\eqref{eq1} and
\(\varepsilon^2=(\alpha R)^{-1}\). It turns out that this problem serves as
a good model for the study of the original Orr--Sommerfeld
problem~\eqref{eq1}--\eqref{eq2}. Now the question is: for which function
\(q(x)\) we can describe the spectral portraits of
problem~\eqref{eq3},~\eqref{eq4} as \(\varepsilon\to 0\)?
The answer is not simple even for particular functions \(q(x)\), for example
\(q(x)=x^2\) or \(q(x)=(x-\alpha)^2\), \(\alpha\in\mathbb R\). Only in the
case \(q(x)=x\) the problem can be solved explicitely~\cite{ci1},~\cite{ci2}.
At the moment we have general result only for analytic monotonous functions.

We shall say that a function \(q(x)\) belongs to the class \(AM\) if the
following conditions hold: \textit{\(q(x)\) is analytic on the segment
\([-1,1]\) and admits an analytic continuation in a domain \(G\) such that the
range \(q(G)\) contains the semistrip \(\Pi=\left\{\lambda\;\vline\;
a<\operatorname{Re}\lambda<b,\;\operatorname{Im}\lambda<0\right\}\) where
\(a=q(-1)\), \(b=q(1)\). Moreover, if \(D=q^{-1}(\Pi)\subset G\), then
\(q(z)\) is continuous on the closure \(\overline D\) and the map
\(q(z):\overline D\to\overline{\Pi}\) is one to one.}

Assuming that \(q(x)\in AM\) we introduce the functions
\[
	Q^{\pm}(\lambda)=\int\limits_{\xi_{\lambda}}^{\pm 1}
	\sqrt{i(q(\xi)-\lambda)}\,d\xi,\qquad \lambda\in\overline{\Pi},
\]
where \(\xi_{\lambda}\in\overline{D}\) is the root of the equation
\(q(\xi)-\lambda=0\), and the function
\[
	Q(\lambda)=\int\limits_{-1}^1\sqrt{i(q(x)-\lambda)}\,dx,
	\qquad\lambda\in\overline{\Pi}.
\]
These functions are analytic in \(\Pi\) and continuous in \(\overline{\Pi}\),
provided \(q(x)\in AM\).

The following theorem is a generalization of the previous author's
result~\cite{ci3}.
\begin{theorem}\label{th1}
Let \(q(x)\in AM\). Then the curves \(\tilde\gamma_{\pm}\) determined in the
semistrip \(\Pi\) by the equations \(\operatorname{Re}Q^{\pm}(\lambda)=0\)
have the only intersection point in \(\Pi\), say \(\lambda_0\). The curve
\(\tilde\gamma_{\infty}\) determined in \(\Pi\) by the equation
\(\operatorname{Re} Q(\lambda)=0\), is a function with respect to the
imaginary axis and intersect the curves \(\tilde\gamma^{\pm}\) in the only
point \(\lambda_0\). Denote by \(\gamma_+\), \(\gamma_-\) and
\(\gamma_{\infty}\) the parts of \(\tilde\gamma_+\), \(\tilde\gamma_-\) and
\(\tilde\gamma_{\infty}\) connecting the intersection point \(\lambda_0\) with
the points \(a\), \(b\) and \(-i\infty\), respectively. Denote also
\[
	\Gamma=\gamma_+\cup\gamma_-\cup\gamma_{\infty}.
\]
Then the spectrum of model problem~\eqref{eq3},~\eqref{eq4} is concentrated
along \(\Gamma\) as \(\varepsilon\to 0\), i.~e. \(\forall\delta>0\) there are
no eigenvalues outside the \mbox{\(\delta\)-neigh}\-bour\-hood of \(\Gamma\),
provided that \(\varepsilon<\varepsilon_0=\varepsilon_0(\delta)\). Moreover,
the eigenvalues of the problem near the curves \(\gamma_+\), \(\gamma_-\) and
\(\gamma_{\infty}\) are determined from the conditions
\begin{align*}
	i\int\limits_{\xi_{\lambda}}^1\sqrt{i(q(\xi)-\lambda)}\,d\xi&=
	\varepsilon\pi(k-1/4),&k&\in\mathbb Z,\\
	i\int\limits_{\xi_{\lambda}}^{-1}\sqrt{i(q(\xi)-\lambda)}\,d\xi&=
	\varepsilon\pi(k-1/4),&k&\in\mathbb Z,\\
	i\int\limits_{-1}^1\sqrt{i(q(\xi)-\lambda)}\,d\xi&=
	\varepsilon\pi k,&k&\in\mathbb Z.
\end{align*}
More precizely, let \(\mu_k^+\), \(\mu_k^-\) and \(\mu_k\) are solutions of
these equations lying on the curves \(\gamma_+\), \(\gamma_-\) and
\(\gamma_{\infty}\) outside \mbox{\(\delta\)-neigh}\-bour\-hoods of the points
\(a\), \(b\) and \(\lambda_0\). Then there is a constant \(C=C(\delta)\) such
that each circle of the radius \(C\varepsilon^2\) centered at the points
\(\mu_k^-\), \(\mu_k^+\) and \(\mu_k\) contains the only eigenvalue. The
eigenvalue counting functions along the curves \(\gamma_+\), \(\gamma_-\) and
\(\gamma_{\infty}\) have representations
\begin{align*}
	N_{\pm}(\lambda)&=\pm\dfrac{1}{i\pi\varepsilon}
	\int\limits_{\xi_{\lambda}}^{\pm 1}\sqrt{i(q(\xi)-\lambda)}\,d\xi+
	O(1),\\
	N(\lambda)&=\dfrac{1}{i\pi\varepsilon}
	\int\limits_{-1}^{1}\sqrt{i(q(\xi)-\lambda)}\,d\xi+O(1),
\end{align*}
where \(O(1)\) is uniformly bounded as \(\varepsilon\to 0\) provided that
\(\lambda\) lies outside some neighbourhoods of the points \(a\), \(b\) and
\(\lambda_0\).
\end{theorem}

Our next results clarify the connection of model problem \eqref{eq3}, 
\eqref{eq4} with original Orr-Sommerfeld problem \eqref{eq1}-\eqref{eq2}.

\begin{theorem}\label{th2}
Let \(q(x) \in AM\). Denote by \(\Delta _m(\lambda, \varepsilon)\) and 
\(\Delta _{OS}(\lambda, \varepsilon)\) the characteristic determinants of 
the model and Orr-Sommerfeld problems, respectively. Then
\[
\Delta _{OS}(\lambda, \varepsilon) = \Delta _{m}(\lambda, \varepsilon)
R(\lambda)(1+O(\varepsilon)),
\]
where \(R(\lambda) \neq 0\) for \(\lambda \in \Pi\) and \(|O(\varepsilon)| 
\leq C\varepsilon\) for \(\lambda \in \Pi\setminus\Gamma_{\delta}\). Here 
\(\Gamma_{\delta}\) is the \(\delta\)-neighboyrhood of the limit 
spectral graph \(\Gamma\) determined in Theorem \ref{th1}, and the 
constant C depends only on \(\delta\).
\end{theorem}

The technique, developed in \cite{ci4} allows to prove the 
similarity of the zero counting functions for a pair of holomorphic 
functions, provided that they have similar growth near the curves 
where the zeros are located (in our case on the boundary of
\(\Gamma_{\delta}\)). Hence, Theorem \ref{th1} leads to the 
following result.

\begin{theorem}\label{th3}
Let \(q(x) \in AM \). Then the limit spectral graph of the 
Orr-Sommerfeld problem coincides with the graph \(\Gamma\)
of the corresponding model problem, i.e. \(\forall \delta >0\)
there are no eigenvalues of problem \eqref{eq1}, \eqref{eq2} 
outside the \( \delta \)-neighbourhood of  \(\Gamma\), provided
that \( R>R_0=R_0(\delta)\). Moreover, the eigenvalue counting
functions \(N_{\pm}(\lambda) \) and \( N_{\infty}(\lambda) \) 
along the curves \(\gamma_{\pm} \) and \( \gamma_{\infty} \) 
for problem \eqref{eq1}, \eqref{eq2} have representations
\begin{align*} 
   N_{\pm}(\lambda)&=\pm \frac{1}{i \pi \varepsilon}
   \int\limits_{\xi_{\lambda}}^{\pm 1}
   \sqrt{i(q(\xi)-\lambda)}\,d\xi\,(1+o(1)) ,\\
   N(\lambda)&=\frac{1}{i \pi \varepsilon}
   \int\limits_{-1}^{1} 
   \sqrt{i(q(\xi)-\lambda)}\,d\xi\,(1+o(1))
\end{align*}
where \( o(1)\to 0 \) as \( \varepsilon \to 0 \) uniformly for all
\( \lambda \) lying outside some neighbourhoods of the points 
\(a, b\) and \( \lambda_0 \).
\end{theorem}

Next our goal is to understand the resolvent behaviour of the 
operators corresponding to problems in question. Denote by 
\(z_k = z_k(x, \varepsilon) \) the normed eigenfunctions of the 
model problems (for symplicity we assume that there is no associated
functions, although they do appear for some discrete values of 
\( \varepsilon _j \to 0 \)). It is known that the system 
\(\{z_n\}_1^{\infty} \) forms a Riesz basis in \(L_2(-1,1)\), 
therefore, the Gramm matrix \(\{(z_k, z_j)\}_{k,j=1}^{\infty} \)
generates a bounded and boundely invertible operator 
\( T \) in \(l_2 \). The number 
\begin{equation}\label{eq5}
C=C(\varepsilon)=\|T\| \cdot \|T^{-1}\|
\end{equation}
characterizes "the quality" of the basis  \(\{z_k\}_1^{\infty} \). It is called
the Riesz constant.

Let \( L = L(\varepsilon)\) be the operator in \( L_2(-1,1) \)
associated with model problem \eqref{eq3}, \eqref{eq4}. It can be 
easily proved that
\[
  \|(L-\lambda)^{-1}\| \leqslant \frac{C(\varepsilon)}
  {\operatorname{dist}(\lambda,\sigma(L))}
\]
where \( C(\varepsilon) \) is defined by \eqref{eq5}.
The question is: how rapidly grows the constant 
\( C(\varepsilon) \)?

\begin{theorem} \label{th5} 
Denote by \( \Omega \) the domain bounded by the curves 
\( \gamma_+, \gamma_- \) and the segment \([a,b]\). Then,
given compact \(K \) in \( \Omega \) there are numbers 
\(\tau > 0 \) and \(\varepsilon _0 > 0 \) such that
\( \|(L-\lambda)^{-1}\| \geqslant e^{\tau/\varepsilon}\)
for all \(\lambda \in K\) and \(\varepsilon<\varepsilon_0\).
In particular, the constant \(C(\varepsilon) \) defined in 
\eqref{eq5} grows exponentially, i.e. there is a number 
\(\eta > 0\) such that \( C(\varepsilon) \geqslant
e^{\eta/\varepsilon} \).
\end{theorem}

It was proved in \cite{ci5} that there is no operator in 
\(L_2(-1,1) \) associated with Orr-Sommerfeld problem.
However, there is such an operator in the Sobolev space
\( {\stackrel{\circ}{W}}\vphantom{W}^1_{2}[-1,1] \).
Denote this operator by \(S\)
and consider the operator \(T\) in \(l_2\) generated by 
the Gramm matrix \(\{(y_k, y_j)\}_{k,j=1}^{\infty}\) 
assuming that \(y_k\) are the normed eigenfunctions of
the operator \(S\) in the space
\( {\stackrel{\circ}{W}}\vphantom{W}^1_{2}[-1,1] \)
and the scalar product is taken in the same space.
                                                  
\begin{theorem}\label{th6}
The resolvent \((S-\lambda)^{-1}\) and the constant \(C(\varepsilon)\) defined
by~\eqref{eq5} admit the same estimates from below as in Theorem~\ref{th5}.
\end{theorem}

The estimates from below for the resolvents \((L-\lambda)^{-1}\) and
\((S-\lambda)^{-1}\) confirm the following fact: for sufficiently small
\(\varepsilon>0\) the
pseudospectra of these operators occupies "almost whole" domain \(\Omega\).
So, the last theorems prove some observations of Trefethen~\cite{ci6} and
Reddy, Schmidt and Henningson~\cite{ci7} on the pseudospectra of the
Orr--Sommerfeld operator. We remark also that a weaker estimate
\(C(\varepsilon)\geqslant\tau\varepsilon^{1/3}\) for the Riesz constant was
obtained by Stepin~\cite{ci8}.

\medskip

\end{document}